\documentclass[12pt, twoside, leqno]{article}

\usepackage{amsmath,amsthm}
\usepackage{amssymb}
\usepackage{amsbsy}

\usepackage{enumitem}

\newtheorem{theorem}{Theorem}[section]
\newtheorem{corollary}[theorem]{Corollary}
\newtheorem{lemma}[theorem]{Lemma}
\newtheorem{proposition}[theorem]{Proposition}

\newtheorem{fact}[theorem]{Fact}

\theoremstyle{definition}
\newtheorem{definition}[theorem]{Definition}

\newtheorem*{xA}{Theorem A}
\newtheorem*{xB}{Theorem B}
\newtheorem*{xC}{Corollary C}
\newtheorem*{xD}{Theorem D}
\newtheorem*{xE}{Corollary E}

\newenvironment{proof*}[1]
  {\renewcommand{\proofname}{#1}%
   \begin{proof}}
  {\end{proof}}

\numberwithin{equation}{section}

\begin{document}

\baselineskip=17pt

\title{Haagerup property and Kazhdan pairs via
ergodic infinite measure preserving  actions 
}

\author{Alexandre I. Danilenko\\
B. Verkin Institute for Low Temperature Physics \& Engineering\\
 Ukrainian National Academy of Sciences\\
47 Nauky Ave.\\
 61164,  Kharkiv, UKRAINE\\
E-mail: alexandre.danilenko@gmail.com
}

\date{}

\maketitle

\renewcommand{\thefootnote}{}



\renewcommand{\thefootnote}{\arabic{footnote}}
\setcounter{footnote}{0}

\begin{abstract}
It is shown that a locally compact second countable group $G$ has the Haagerup property if and only if
there exists a sharply weak mixing 0-type  measure preserving free  $G$-action $T=(T_g)_{g\in G}$ on an infinite $\sigma$-finite standard measure space  $(X,\mu)$ admitting an exhausting  $T$-F{\o}lner sequence  (i.e. a sequence $(A_n)_{n=1}^\infty$ of measured subsets of finite measure such that   $A_1\subset A_2\subset\cdots$, $\bigcup_{n=1}^\infty A_n=X$
and $\lim_{n\to\infty}\sup_{g\in K}\frac{\mu(T_gA_n\triangle A_n)}{\mu(A_n)}= 0$ for each compact $K\subset G$).
It is also shown that 
a pair of groups $H\subset G$ has property (T) if and only if there is a $\mu$-preserving  $G$-action $S$ on  $X$ admitting an $S$-F{\o}lner sequence  and such that
$S\restriction H$ is weakly mixing.
These refine some recent results by Delabie-Jolissaint-Zumbrunnen and Jolissaint.
\end{abstract}

\section{Introduction}

Throughout this paper $G$ is a non-compact locally compact second countable group.
It  has {\it the Haagerup property} if there is  a weakly continuous unitary representation $V$ of $G$ in a separable Hilbert space $\mathcal H$ such that $\lim_{g\to\infty }V(g)= 0$  in the weak operator topology and

\begin{enumerate}[label=($\ast$)]
\item \label{ast}
for each $\epsilon>0$ and every compact subset $K\subset G$, there is a unit vector $\xi\in\mathcal H$
 such that $\sup_{g\in K}\|V(g)\xi-\xi\|<\epsilon$. 
\end{enumerate}
 The amenable groups, $SO(n,1)$ and $SU(n,1)$ for  each $n\ge 2$, the free groups, the Coxeter groups have the  Haagerup property \cite{Ch--Va}.
The class of discrete  countable Haagerup groups  is closed under free products and wreath products \cite{CoStVa}.
For more information about the Haagerup property we refer to \cite{Ch--Va}.
There is a purely dynamical description  of this property: $G$ is  Haagerup  if and only if there exists a mixing  non-strongly ergodic probability preserving free
$G$-action  \cite[Theorem~2.2.2]{Ch--Va} (see \S\ref{sec2} for the definitions).
Recently, an infinite measure preserving counterpart of this result was discovered  in \cite{DeJoZu}:

\begin{xA} $G$ has the Haagerup property if and only if there is a 0-type measure preserving  $G$-action $T=(T_g)_{g\in G}$ on an infinite $\sigma$-finite measure space $(X,\frak B,\mu)$ admitting a sequence of non-negative unit vectors $(\xi_n)_{n=1}^\infty$ in $L^2(X,\mu)$ such that $\lim_{n\to\infty}\sup_{g\in K}\langle\xi_n\circ T_g,\xi_n\rangle=1$ for each compact $K\subset G$.
\end{xA}

We recall that $T$ is called {\it of 0-type} if $\lim_{g\to\infty}\mu(T_gA\cap B)=0$ for all subsets $A,B\in\frak B$ of finite measure.
In this paper we provide a much shorter alternative proof of Theorem~A which is grounded on the Moore-Hill concept of restricted infinite products of probability measures \cite{Hi}.

 We note that the  0-type  for infinite measure preserving systems is a natural counterpart of the mixing 
for probability preserving systems.
However  unlike mixing,  the 0-type is not a ``strong'' asymptotic property.
It  implies neither weak mixing nor ergodicity.
Moreover, the  totally dissipative actions are  all of 0-type. 
In view of that  the description in Theorem~A does not look sharp  from the ergodic theory point of view.
Our first main result in this work is the following 
 finer  ergodic criterion of the Haagerup property.\footnote{Conservativeness, ergodicity, weak mixing and sharp weak mixing are not spectral invariants of the underlying dynamical systems. Hence the principal difference of Theorem B from Theorem~A is that it provides non-spectral ergodic characterization  of the Haagerup property.}

\begin{xB} The following are equivalent.
\begin{enumerate}[label=\upshape(\roman*), leftmargin=*, widest=iii]
 \item
 $G$ has the Haagerup property.
 \item
 There exists a sharply weak mixing (conservative) 0-type measure preserving free $G$-action $T$ on an infinite $\sigma$-finite standard measure space admitting an exhausting  $T$-F{\o}lner sequence of subsets.
\item
 There exists a  sharply weak mixing (conservative) 0-type  measure preserving free  $G$-action $T$ on an infinite $\sigma$-finite standard measure space $(X,\frak B,\mu)$ admitting a  $T$-F{\o}lner sequence $(A_n)_{n=1}^\infty$ such that $\mu(A_n)=1$ for all $n\in\Bbb N$.
\end{enumerate}
 \end{xB}

We say that $(A_n)_{n=1}^\infty$ is {\it $T$-F{\o}lner} if $\mu(A_n)<\infty$ and
    $$
 \sup_{g\in K}\frac{\mu(A_n\triangle T_gA_n)}{\mu(A_n)}\to0\qquad\text{as $n\to\infty$}
       $$
        for each compact subset $K\subset G$.
       If $A_1\subset A_2\subset\cdots$ and $\bigcup_{n=1}^\infty A_n=X$, we say that $(A_n)_{n=1}^\infty$ is {\it exhausting}.
We note that sharp weak mixing (see \S\ref{sec2} for the definition) implies ergodicity and weak mixing.
To prove (the non-trivial part of) Theorem~B, we apply the Moore-Hill construction \cite{Hi} to the mixing non-strongly ergodic $G$-action from \cite[Theorem~2.2.2]{Ch--Va}  (cf. the construction of  $II_\infty$ ergodic Poisson suspensions of countable amenable groups from \cite{DaKo}).
Then we observe that the action $T$ that we obtain is IDPFT (see \S\ref{sec2} and \cite{DaLe}, where such actions  were introduced). 
Hence, by the properties of IDPFT systems, $T$ is sharply weak mixing whenever we show that it is conservative. 
To show the conservativeness of $T$ it remains to choose the parameters of the Moore-Hill construction in an a appropriate way.

As a corollary from Theorem~B, we obtain one more dynamical characterization of the Haagerup property in terms of  Poisson actions.

\begin{xC}  $G$ has the Haagerup property if and only if there exists a mixing (probability preserving) Poisson $G$-action that is not strongly ergodic.
 \end{xC}

 Our next purpose is to obtain a  ``parallel'' characterization  of property~(T) which is a reciprocal to the Haagerup property.
 We recall \cite[ Definition 1.1]{Jo1} that given a non-compact closed subgroup $H$ of $G$, the pair $H\subset G$ {\it has property (T)} if for each unitary representation $V$ of $G$ satisfying~\ref{ast},
 there is a unit vector which is invariant under $V(h)$ for every $h\in H$.
 The property (T) for a single group $G$ corresponds to the (T)-property for the pair $G\subset G$.
 Using the techniques developed for proving~Theorem~B we obtain  an ergodic  (non-spectral)  characterization of Kazhdan pairs  that refines  a spectral characterization from  \cite{Jo2}.

 \begin{xD} 
\begin{enumerate}[label=\upshape(\roman*), leftmargin=*, widest=ii]
 \item
 Assume that a  pair $H\subset G$ has property (T).
Let $S=(S_g)_{g\in G}$   be a  measure preserving $G$-action 
 on a $\sigma$-finite infinite standard measure space $(Y,\frak C,\nu)$, such that $S\restriction H:=(S_h)_{h\in H}$ has no invariant subsets of  positive finite measure.
  Then this action admits no $S$-F{\o}lner sequences.
  \item
 If a  pair $H\subset G$ does not have property (T) then there is a measure preserving $G$-action $S$
 on a $\sigma$-finite infinite measure space which has an exhausting  
$ S$-F{\o}lner sequence
 and
  such that  
$ S\restriction H$ is  weakly mixing.
 \end{enumerate}
 \end{xD}

Let us say that  $S\restriction H$ is {\it of  weak 0-type} if  there is a subsequence $h_n\to\infty$ in $H$ such that $\lim_{n\to\infty}\nu(S_{h_n}A\cap B)= 0$ for all subsets $A,B\in\frak C$ of finite measure.
 Then replacing   the ``has no invariant subsets of positive finite measure" in (i) with  a stronger ``is of  weak 0-type", 
 and
  the ``weakly mixing'' in (ii) with  a weaker ``of  weak 0-type''  we obtain exactly \cite[Theorem~1.5]{Jo2}.
 
 \begin{xE}
 A  pair $H\subset G$ has property (T) if and only if 
every (probability preserving)  Poisson $G$-action   
with weakly mixing  $H$-subaction is strongly ergodic.
The same is also true with ``ergodic'' in place of ``weakly mixing".
  \end{xE}

 The outline of the paper is as follows.
 In Section~\ref{sec2} we state all necessary definitions related to the basic dynamical concepts of group actions both  in  the nonsingular and   and  finite measure preserving cases, restricted infinite powers of probability measures, IDPFT actions and Poisson actions.
 In Section~\ref{sec3} we prove Theorems~B and Corollary~C.
 Section~\ref{sec4} is devoted  to the proof of Theorems~D and Corollary~E.

\section{Definitions and preliminaries}\label{sec2}

\subsection*{Nonsingular and measure preserving $G$-actions}
Nonsingular actions appear in the proof of Theorem~B.
We remind several basic concepts related to them (see \cite{Aa}, \cite{ScWa}, \cite{GlWe2}, \cite{DaKo}).

\begin{definition}\label{def2.1} Let $S=(S_g)_{g\in G}$ be a nonsingular $G$-action on a standard  probability space $(Z,\frak F,\kappa)$.
\begin{enumerate}[label=\upshape(\roman*), leftmargin=*, widest=iii]
 \item
$S$ is called {\it totally dissipative} if the partition of $Z$ into the $S$-orbits is measurable and the $S$-stabilizer
of a.e. point is compact, i.e. there is a measurable subset of $Z$ which meets a.e. $S$-orbit exactly once, and for a.e. $z\in Z$, the subgroup $\{g\in G\mid S_gz=z\}$ is compact in $G$.
\item
$S$ is called {\it conservative} if there is no  $S$-invariant subset $A\subset Z$ of positive measure such that the restriction of $S$ to $A$ is totally dissipative.
\item
There is a unique (mod 0) partition of $X$ into two invariant subsets $\mathcal D(S)$ and $\mathcal C(S)$ such that $S\restriction\mathcal D(S)$ is totally dissipative and $S\restriction\mathcal D(S)$ is conservative. We call $\mathcal D(S)$ and $\mathcal C(S)$ the {\it dissipative} and {\it conservative part} of $S$ respectively.
\item
$S$ is called {\it ergodic} if each measurable $S$-invariant subset of $Z$ is either $\mu$-null or $\mu$-conull.
\item $S$ is called {\it weakly mixing} if for each ergodic probability preserving 
$G$-action $R=(R_g)_{g\in G}$, the product $G$-action $(S_g\times R_g)_{g\in G}$ is 
ergodic.
\item
$S$ is called {\it sharply weak mixing} \cite{DaLe} if $S$ is  conservative, ergodic and for each  ergodic 
conservative nonsingular $G$-action $R=(R_g)_{g\in G}$, 
 the product $G$-action $(S_g\times R_g)_{g\in G}$ is either  ergodic or  totally dissipative.
 \end{enumerate}
\end{definition}

We also remind some concepts related to finite measure preserving actions.

\begin{definition}\label{def 2.2} Suppose that $\kappa(Z)=1$ and $\kappa\circ S_g=\kappa$
for all $g\in G$.
\begin{enumerate}[label=\upshape(\roman*), leftmargin=*, widest=iii]
   \item 
   $S$ is called {\it mixing} if $\lim_{g\to\infty}\kappa(S_gA\cap B)=\kappa(A)\kappa(B)$
   for all $A,B\in\frak F$.
        \item
        A sequence of Borel subsets $(A_n)_{n=1}^\infty$ in $X$  of strictly positive measure is called {\it $T$-asymptotically invariant} if  for each compact subset $K\subset G$, we have that $\sup_{g\in K}\kappa(A_n\triangle T_gA_n)\to 0$ as $n\to\infty$.
     \item
      $T$ is called {\it strongly ergodic} if each $T$-asymptotically invariant sequence $(A_n)_{n=1}^\infty$ is trivial, i.e. $\lim_{n\to\infty}\kappa(A_n)(1-\kappa(A_n))=0$. 
  \end{enumerate}
  \end{definition}

 We now state a corollary from the Schmidt-Walters theorem  \cite[Theorem~2.3]{ScWa}
 as it appeared in \cite[Theorem~7.3]{ArIsMa}.
 For a detailed proof of a sharper result, we refer to  \cite[Theorem~7.14]{ArIsMa}. 
 
 \begin{lemma}\label{lem2.3}
  Let $S=(S_g)_{g\in G}$ be a mixing measure preserving action on a standard probability space $(Y,\frak C,\nu)$ and
 let $R=(R_g)_{g\in G}$
be a   conservative nonsingular $G$-action on a standard probability space $(Y,\frak C,\nu)$.
If 
$$
F:Y\times Z\to\Bbb C
$$
 is an $(S\times R)$-invariant Borel function then there exists a Borel $R$-invariant function $f:Z\to\Bbb C$ such that
$F(y,z)=f(z)$  at a.e. $(y,z)\in Y\times Z$.
 \end{lemma}

 \begin{corollary}\label{cor2.4}  Let $S=(S_g)_{g\in G}$ be a mixing measure preserving action on a standard probability space $(Y,\frak C,\nu)$ and let $R=(R_g)_{g\in G}$
be a   nonsingular $G$-action on a standard probability space $(Z,\frak D,\kappa)$. 
Then 
$$\mathcal D(S\times R)=Y\times \mathcal D(R)\text{ and  }
\mathcal C(S\times R)=Y\times \mathcal C(R).
$$
 \end{corollary}

\subsection*{Restricted infinite powers of probability measures}
Let $(Y,\frak C,\gamma)$ be a  standard non-atomic probability space. 
Fix a sequence  $\boldsymbol B:=(B_n)_{n=1}^\infty$ of subsets from $\frak C$ of  positive measure.
Let $(X,\frak B):=(Y,\frak C)^{\otimes\Bbb N}$.
For each $n\in\Bbb N$, we set
 $
 \boldsymbol B^n:=Y^n\times B_{n+1}\times B_{n+2}\times\cdots\in\frak B.
 $
 Then 
 $$
  \boldsymbol B^1\subset \boldsymbol B^2\subset\cdots.
  $$
We define a measure $\gamma^{ \boldsymbol B}$ on  $(X,\frak B)$ by 
the following sequence of restrictions (see  \cite{Hi} for details):
 $$
 \gamma^{ \boldsymbol B}\restriction\boldsymbol B^{n}:=\frac{\gamma}{\gamma(B_1)}\otimes
 \cdots\otimes \frac{\gamma}{\gamma(B_n)}\otimes
 \frac{\gamma\restriction B_{n+1}}{\gamma(B_{n+1})}\otimes
  \frac{\gamma\restriction B_{n+2}}{\gamma(B_{n+2})}\otimes\cdots.
 $$
 Since the restrictions are compatible, $\gamma^{ \boldsymbol B}$ is well defined.
 We note that $\gamma^{ \boldsymbol B}$ is supported on the subset 
 $\bigcup_{n=1}^\infty \boldsymbol B^n\subset Y$ 
 and  $\gamma^{ \boldsymbol B}(\boldsymbol B^n)=\prod_{j=1}^n\gamma(B_j)^{-1}$ for each $n$.
Hence, $\gamma^{ \boldsymbol B}$ is  $\sigma$-finite.
It is infinite if and only if $\prod_{n=1}^\infty \gamma(B_n)=0$. 
 
 \begin{definition}\label{def2.5} We call $\gamma^{ \boldsymbol B}$  {\it the restricted infinite power of $\gamma$ with respect to $\boldsymbol B$.}
 \end{definition}

Given a   $\gamma$-preserving Borel bijection $T$ of $Y$, we let
$\boldsymbol T:=\bigotimes_{n=1}^\infty T$ and 
$\boldsymbol T\boldsymbol B:=(TB_n)_{n=1}^\infty$.
A straightforward verification shows that $\gamma^{ \boldsymbol B}\circ T^{-1}=\gamma^{\boldsymbol T \boldsymbol B}$.

 \begin{proposition}\label{pro2.6} If $\sum_{n=1}^\infty\frac{\gamma(B_n\triangle TB_n)}{\gamma(B_n)}<\infty$ then
 $\boldsymbol T$ preserves $\gamma^{ \boldsymbol B}$.
 \end{proposition}
 
 \begin{proof} For each $n\in\Bbb N$ and an arbitrary  Borel subset $A\subset Y^n$, we let 
 $A':=A\times B_{n+1}\times B_{n+2}\times\cdots\in\frak B$.
 Then for every $m>n$,
 $$
 (\gamma^{\boldsymbol T \boldsymbol B}\restriction(\boldsymbol T\boldsymbol B)^{m})(A')=\frac{\gamma^{\otimes n}(A)}{\gamma(B_1)\cdots\gamma(B_n)}\prod_{j>m}\frac{\gamma(B_j\cap TB_j)}{\gamma(B_j)}.
 $$
 Passing to the limit as $m\to\infty$, we obtain that 
 $$
 \gamma^{ \boldsymbol T\boldsymbol B}(A')=\lim_{m\to\infty}
 (\gamma^{\boldsymbol T \boldsymbol B}\restriction(\boldsymbol T\boldsymbol B)^{m})(A')=\frac{\gamma^{\otimes n}(A)}{\gamma(B_1)\cdots\gamma(B_n)}=\gamma^{ \boldsymbol B}(A').
 $$
 Hence $\gamma^{ \boldsymbol B}\circ T^{-1}= \gamma^{ \boldsymbol T\boldsymbol B}=\gamma^{ \boldsymbol B}$, as desired.
 \end{proof}
 
 We note that under the condition of Proposition~\ref{pro2.6}, 
 $$
  \gamma^{ \boldsymbol B}(\boldsymbol T\boldsymbol B^n\cap\boldsymbol B^n)={\prod_{j=1}^n\gamma(B_j)^{-1}}
 \prod_{j>n}\frac{\gamma(B_j\cap TB_j)}{\gamma(B_j)}.
 $$
 Hence 
 \begin{equation}\label{2-1}
\lim_{n\to \infty }\frac{\gamma^{ \boldsymbol B}(\boldsymbol T\boldsymbol B^n\cap\boldsymbol B^n)}{\gamma^{ \boldsymbol B}(\boldsymbol B^n)}= 1.
 \end{equation}
    
  We note that  the Hilbert space $L^2(X,\gamma^{\boldsymbol B})$  is  the infinite tensor product of the sequence of Hilbert spaces $(L^2(Y,\gamma(B_n)^{-1}\gamma))_{n=1}^\infty$
 along the stabilizing sequence $(1_{B_n})_{n=1}^\infty$ of unit vectors (see \cite{Gu} for the definition).
    In particular, we see that the linear subspace
    $$
    \bigcup_{n=1}^\infty\Bigg\{f\otimes \bigotimes_{j>n}1_{B_j}\,\bigg|\, f\in L^2(Y^n,\gamma^{\otimes n})\Bigg\}
    $$
is dense in  $L^2(X,\gamma^{\boldsymbol B})$.
Let $U_T$ and $U_{\boldsymbol T}$ stand for the Koopman operators associated to  $T$ and 
$\boldsymbol T$ in $L^2(Y,\gamma)$ and $L^2(X,\gamma^{\boldsymbol B})$ respectively.
The following proposition is verified straightforwardly.

\begin{lemma}\label{lem2.7}
For each $n\in\Bbb N$ and arbitrary functions $f,r\in L^2(Y^n,\gamma^{\otimes n})$,
we have that
 \begin{gather*}
 U_{\boldsymbol T} \bigg(f\otimes \bigotimes_{j>n}1_{B_{j}}\bigg)=\lim_{m\to\infty}(U_{T})^{\otimes n}f\otimes \bigg(\bigotimes_{j=n+1}^m1_{TB_{j}}\bigg)\otimes 
 \bigotimes_{j>m}1_{B_j}\quad\text{ and}\\
\bigg\langle U_{\boldsymbol T}\bigg(f\otimes \bigotimes_{j>n}1_{B_{j}}\bigg),r\otimes \bigotimes_{j>n}1_{B_{j}}\bigg\rangle=\frac{\langle(U_{T})^{\otimes n}f,r\rangle}{\prod_{j=1}^n\gamma(B_j)}\prod_{j>n}\frac{\gamma(TB_j\cap B_j)}{\gamma(B_j)}.
\end{gather*}
The above limit is considered in the strong topology in $L^2(X,\gamma^{\boldsymbol B})$.
\end{lemma}

\subsection*{IDPFT-actions}
IDPFT-actions were introduced in \cite{DaLe} in the case, where $G=\Bbb Z$.
In \cite{DaKo}, IDPFT actions of arbitrary discrete countable groups were studied.
In this paper we consider IDPFT-actions for  arbitrary locally compact second countable groups.

\begin{definition}\label{def2.8} Let $T_n=(T_n(g))_{g\in G}$ be an ergodic  measure preserving $G$-action
on a standard probability space $(Y_n,\frak C_n,\nu_n)$, let $\mu_n$ be a probability measure on $\frak C_n$ and let $\mu_n\sim\nu_n$ for each $n\in\Bbb N$.
We put $(X,\frak B,\mu):=\bigotimes_{n=1}^\infty(Y_n,\frak C_n,\mu_n)$,  $T(g):=\bigotimes_{n=1}^\infty T_n(g)$  and $T:=(T(g))_{g\in G}$.
If $\mu\circ T(g)\sim\mu$ for each $g\in G$ then the nonsingular dynamical system $(X,\frak B,\mu, T)$
is called  an {\it infinite direct product of finite types (IDPFT)}.
\end{definition}

We will need the following fact, extending \cite[Proposition~2.3]{DaLe} from $\Bbb Z$-actions
to arbitrary $G$-actions.

\begin{proposition}\label{pro2.9} Let $(X,\frak B,\mu, T)$ be an IDPFT system as in Definition~\ref{def2.8}.
If $T_n$ is mixing for each $n\in\Bbb N$ then $T$ is either  totally dissipative  or conservative.
If $T$ is conservative then
$T$ is sharply weak mixing.
\end{proposition}

\begin{proof} It follows from  Corollary~\ref{cor2.4}  that 
$$
\mathcal D(T)=Y_1\times\cdots \times Y_n\times\mathcal D\bigg(\bigotimes_{j>n} T_j\bigg) \quad\mod 0
$$
for each $n>0$.
By Kolmogorov's 0-1 law that either $\mu(\mathcal D(T))=0$ and hence $T$ is conservative or $\mu(\mathcal D(T))=1$ and hence $T$ is totally dissipative.
Thus, the first claim is proved.
We do not provide a  proof for the second claim because it is an almost verbal repetition of the proof of \cite[Proposition~2.3]{DaLe}:  just replace the reference to \cite[Theorem~B]{DaLe} there with a reference to Lemma~\ref{lem2.3}.
Of course, $\mu$ is not concentrated on a single orbit.
\end{proof}

\subsection*{Poisson suspensions \rm{(see \cite{LaPe} and \cite{Ro} for details)}}
Let $(X,\frak B)$ be a standard Borel space and let $\mu$ be an infinite $\sigma$-finite 
non-atomic measure on $X$.
 Let $X^*$ be the set of  $\Bbb Z_+$-valued  ($\sigma$-finite) measures on $X$.
  For each subset $A\in\frak B$ with $0<\mu(A)<\infty$, we define a mapping $N_A:X^*\to\Bbb R\cup\{+\infty\}$ by setting
 $N_A(\omega):=\omega(A)$.
 Let $\frak B^*$ stand for the smallest $\sigma$-algebra on $X^*$ such that the mappings
 $N_A$ are all $\frak B^*$-measurable.
There is a unique probability measure $\mu^*$  on $(X^*,\frak B^*)$ satisfying the following two conditions:
\begin{itemize}
\item 
 the measure $\mu^*\circ N_A^{-1}$ is the Poisson distribution
with parameter $\mu(A)$ for each $A\in\frak B$ with $\infty>\mu(A)>0$,
\item
 given a finite family $A_1,\dots, A_q$ of mutually disjoint subsets $A_1,\dots, A_q$ of $X$ of
finite positive measure, the corresponding random variables $N_{A_1},\dots, N_{A_q}$ defined on the space $(X^*,\frak B^*, \mu^*)$ are independent.
\end{itemize}
Then $(X^*,\frak B^*, \mu^*)$ is a Lebesgue space.
The mapping $N_A$ is finite $\mu^*$-almost everywhere for each $A\in\frak B$ with $\infty>\mu(A)>0$.
Moreover, for $\mu^*$-a.e. $\omega$, there exist countably many points $x_j\in X$, $j\in\Bbb N$, such that $\omega=\sum_{j\in\Bbb N}\delta_{x_j}$. 
For each $\mu$-preserving $G$-action $T=(T_g)_{g\in G}$, we define a $G$-action $T^*=(T_g^*)_{g\in G}$  on $(X^*,\frak B^*, \mu^*)$ by setting 
$$
T_g^*\omega:=\omega\circ T_g^{-1} \qquad\text{for all $\omega\in X^*$ and $g\in G$.}
$$
Then $T^*$ preserves $\mu^*$.

\begin{definition}\label{def2.10}
The dynamical     system $(X^*,\frak B^*, \mu^*, T^*)$ is called
 {\it the Poisson suspension} of $(X,\frak B, \mu, T)$.
 A probability preserving $G$-action
  is called {\it Poisson} if it  is isomorphic to a Poisson suspension of some infinite $\sigma$-finite measure preserving $G$-action.
  \end{definition}
 
\subsection*{Unitary representations of $G$ and Koopman representations of measure preserving actions}
Let $V=(V(g))_{g\in G}$ be a weakly continuous unitary representation of $G$ in a separable Hilbert space $\mathcal H$. 
We will always assume that $V$ is a complexification of an orthogonal representation of $G$ in a real Hilbert space.

\begin{definition}\label{def2.11} $V$ is called {\it weakly mixing} if $V$ has no nontrivial finite dimensional invariant subspaces.
\end{definition}

  The {\it Fock space $\mathcal F(\mathcal H)$ over $\mathcal H$} is the orthogonal sum $\bigoplus_{n=0}^\infty\mathcal H^{\odot n}$, where $\mathcal H^{\odot n}$ is the $n$-th symmetric tensor power of $\mathcal H$ when $n>0$ and 
 $\mathcal H^{\odot 0}:=\Bbb C$.
 By $\exp V=(\exp V(g))_{g\in G}$ we denote the corresponding unitary representation of $G$ in $\mathcal F(\mathcal H)$, i.e.
 $\exp V(g):=\bigoplus_{n=0}^\infty V(g)^{\odot n}$ for each  $g\in G$ (see \cite{Gu} for details).

 Let $T=(T_g)_{g\in G}$ be a measure preserving $G$-action  on a $\sigma$-finite nonatomic standard measure space $(X,\frak B,\mu)$.
 Denote by $U_T=(U_T(g))_{g\in G}$ the associated (weakly continuous) unitary {\it Koopman representation} of $G$
 in $L^2(X,\mu)$:
 $$
 U_T(g)f:=f\circ T_g^{-1},\quad\text{for all }g\in G.
 $$
 If $\mu(X)<\infty$, we let 
 $$
  L^2_0(X,\mu):= L^2(X,\mu)\ominus\Bbb C=\left\{f\in L^2(X,\mu)\mid\int_Xfd\mu=0\right\}.
  $$
  We will need the following fact.
  
  \begin{fact}\label{fac2.12} Let $V$ and $(X,\frak B,\mu,T)$ be as above in this subsesction.
  \begin{enumerate}[label=\upshape(\roman*), leftmargin=*, widest=iii]
   \item 
 If $\mu(X)<\infty$ then $T$ is weakly mixing if and only if
  $U_T\restriction L^2_0(X,\mu)$ is weakly mixing {\rm \cite{BeRo}.}
     \item
     $V$ is weakly mixing if and only if $(\exp V)\restriction(\mathcal F(\mathcal H)\ominus \Bbb C)$
 is weakly mixing {\rm \cite[Theorem~A3]{GlWe1}}.
  \item
  If $\mu(X)=\infty$ then $U_{T^*}$ is canonically unitarily equivalent to $\exp U_T$ {\rm \cite{Ro}}.
 \item
If $(Y,\frak C,\nu, S)$ denote the Gaussian dynamical system ($G$-action) associated with $V$ then
 $U_S$ is canonically unitarily equivalent to $\exp V$ {\rm \cite{Gu}}.
   \item
 $V$ is weakly mixing if and only if there is a sequence $g_n\to\infty$ in $G$
   such that $V(g_n)\to 0$ weakly as $n\to\infty$ {\rm \cite[ Corollary~1.6, Theorem~1.9]{BeRo}}.
  \end{enumerate}
  \end{fact}

\section{The Haagerup property}\label{sec3}
 
  In this section we prove Theorem~B (and hence Theorem~A) and Corollary~C.
 Prior to this we state  a folklore proposition.
 
  \begin{proposition}\label{pro3.1} Let $S=(S_g)_{g\in G}$ be a measure preserving $G$-action  on an infinite   $\sigma$-finite standard measure space $(Z,\frak Z,\kappa)$.
  Let $\frak Z_0\subset \frak Z$ stand for the ring of subsets of finite measure.
Let  $Z_1\subset Z_2\subset\cdots$ be a sequence of subsets from $\frak Z_0$ such that
\begin{itemize}
\item
$\bigcup_{n=1}^\infty Z_n=Z$ and
\item
for each $n>0$, there is a finite partition $\mathcal P_n$ of $Z_n$
into subsets of equal measure such that the sequence $(\mathcal P_n)_{n=1}^\infty$ approximates $(\frak Z,\kappa)$,
i.e.
for each  $B\in \frak Z_0$ and $\epsilon>0$, there is $N>0$ such that  if $n>N$ then there exists  a $\mathcal P_n$-measurable subset $B_n$ with $\kappa(B\triangle B_n)<\epsilon$.
\end{itemize}
\begin{enumerate}[label=\upshape(\roman*), leftmargin=*, widest=iii]
\item
If  for each compact subset $K\subset G$, an integer $n>0$ and a $\mathcal P_n$-atom $A$, there exist a finite family
  $g_1,\dots,g_l\in G\setminus K$ and  mutually disjoint subsets $A_1,\dots,A_l$ of $A$ such that $\bigsqcup_{i=1}^lS_{g_i}A_i\subset A$    and $\kappa(\bigsqcup_{i=1}A_i)>0.1\kappa(A)$
 then   $S$ is conservative. 
  \item
 If   for each $n>0$ and every pair of $\mathcal P_n$-atoms $A$ and $B$,
  there exist a finite family
  $g_1,\dots,g_l\in G$ and  mutually disjoint subsets $A_1,\dots,A_l$ of $A$ such that $\bigsqcup_{i=1}^lS_{g_i}A_i\subset B$    and $\kappa(\bigsqcup_{i=1}A_i)>0.1\kappa(A)$ then
    $S$ is ergodic.
     \end{enumerate}
  \end{proposition}

  \begin{proof*}{Idea of the proof}
  (i) Take a subset $A\subset Z$ of positive finite  measure and a compact subset $K\subset G$.
  Our purpose is to find $g\in G\setminus K$ such that $\mu(A\cap S_gA)>0$.
  For that, we select $n>0$ and a $\mathcal P_n$-measurable subset $B_n$ such that $\mu(A\triangle B_n)\le 0.001\mu(A)$. 
  Then there is a $\mathcal P_n$-atom $B\subset B_n$ such that $\mu(A\cap B)>0.999\mu(B)$.
  Utilizing the condition of (i), we can find a subset $C\subset B$ and $g\in G\setminus K$
  such that $\mu(A\cap C)>0.9\mu(C)$ and $\mu(A\cap S_gC)>0.9\mu(C)$.
  This yields that $\mu(A\cap S_gA)>\mu((A\cap C)\cap (S_gA\cap C))>0$.
  
  (ii) is proved in a similar way.
  \end{proof*}

   \begin{proof*}{Proof of Theorem~B} 
   The implications  (ii)$\Rightarrow$(i) and  (iii)$\Rightarrow$(i) are trivial.

 We now prove (i)$\Rightarrow$(ii).
 Let $G$ have the Haagerup property.
 By \cite[Theorem~2.2.2]{Ch--Va}, there is a mixing measure preserving free $G$-action $T=(T_g)_{g\in G}$ on a standard probability space $(Y,\frak C,\gamma)$ and a $T$-asymptotically invariant sequence $\boldsymbol B:=(B_n)_{n=1}^\infty$ such that $\gamma(B_n)=0.5$ for each $n\in\Bbb N$.
 Passing to a subsequence, if necessary, we may (and will) assume that
 for each compact subset $K\subset G$, 
 \begin{equation}\label{3-1}
 \sum_{n=1}^\infty\sup_{g\in K}\gamma(T_gB_n\triangle B_n)<+\infty.
 \end{equation}
We now let $(X,\frak B):=(Y,\frak C)^{\otimes\Bbb N}$.
Endow this standard Borel space with the  ($\sigma$-finite) restricted infinite power 
$\gamma^{\boldsymbol B}$ of $\gamma$ with respect to $\boldsymbol B$.
 Since $\prod_{n\in\Bbb N}\gamma(B_n)=0$, it follows that $\gamma^{\boldsymbol B}(X)=\infty$. 
For $g\in G$, let $\boldsymbol T_g:=T_g^{\otimes \Bbb N}$. 
In view of (\ref{3-1}), it follows from Proposition~\ref{pro2.6} that  $\gamma^{\boldsymbol B}\circ \boldsymbol T_g=\gamma^{\boldsymbol B}$.
Thus,
$\boldsymbol T:=(\boldsymbol T_g)_{g\in G}$ is a measure preserving $G$-action on $(X,\frak B,\gamma^{\boldsymbol B})$.
Since the map $X\ni x=(y_n)_{n=1}^\infty\to y_1\in Y$ intertwines $\boldsymbol T$ with $T$ and $T$ is free, $\boldsymbol T$ is free too.
For each $n\in\Bbb N$, define a subset $\boldsymbol B^n\subset X$ in the same way as in~\S1.
It follows from  (\ref{2-1}) that the sequence
$(\boldsymbol B^n)_{n\in \Bbb N}$ is $\boldsymbol T$-F{\o}lner.
Moreover, it is exhausting and $\gamma^{ \boldsymbol B}(\boldsymbol B^n)=2^n$ for each $n$.
To show that $\boldsymbol T$ is of 0-type,
we first note that  since $T$ is mixing then for each pair of integers $n<m$,
$$
\lim_{g\to\infty}\prod_{j=n}^{m}\frac{\gamma(T_gB_j\cap B_j)}{\gamma(B_j)}=\prod_{j=n}^m\gamma(B_j)= 2^{-m+n-1}.
$$
Hence given two functions $f,r\in L^2(Y^n,\gamma^{\otimes n})$, we deduce from Lemma~\ref{lem2.7} that
$$
\lim_{g\to\infty}\bigg\langle U_{\boldsymbol T}\bigg(f\otimes \bigotimes_{j>n}1_{B_{j}}\bigg),r\otimes \bigotimes_{j>n}1_{B_{j}}\bigg\rangle= 0.
$$
This implies that $\boldsymbol T$ is of 0-type.\footnote{It is worthy to note that at this point we have proved completely Theorem~A.}

Thus, it remains to show that $\boldsymbol T$ is sharply weak mixing.
To this purpose, we first prove
that  upon a replacement of $\gamma^{\boldsymbol B}$ with an equivalent probability measure,   $\boldsymbol T$ is an IDPFT. 
 For that, we choose
 a sequence of  reals $(\epsilon_n)_{n=1}^\infty$ such that $0<\epsilon_n<1$ for each $n$ and $\sum_{n=1}^\infty\epsilon_n<\infty$.
  Then we define, for each $n\in\Bbb N$, a Borel function $\phi_n:Y\to\Bbb R_+$ by setting 
  $$
  \phi_n:=2\epsilon_n1_{X\setminus B_n}+2(1-\epsilon_n)1_{B_n}.
  $$
   Since $\gamma(B_n)=\frac12$, a simple verification yields  that  $\int_Y\phi_nd\gamma=1$.
Denote by   $\mu_n$  the probability measure on $Y$ such that  $\mu_n\sim\gamma$ and 
   $\frac{d\mu_n}{d\gamma}:=\phi_n$.
  We now recall that 
  given two probability measures $\alpha$ and  $\beta$ on $(Y,\frak C)$ such that $\alpha\prec\gamma$ and $\beta\prec\gamma$, the squared {\it Hellinger distance} between $\alpha$ and $\beta$ is
  $$
  H^2(\alpha,\beta) :=\frac 12\int_Y\left(\sqrt{\frac{d\alpha}{d\gamma}}-\sqrt{\frac{d\beta}{d\gamma}}\right)^2d\gamma.
  $$
A straightforward computation yields that
  $$
  H^2\left(\frac1{\gamma(B_n)}\gamma\restriction B_n,\mu_n\right)=\frac 12\int_Y(\sqrt 2\cdot 1_{B_n}-\sqrt{\phi_n})^2d\gamma=\frac 1{2}((\sqrt{1-\epsilon_n}-1)^2+\epsilon_n)
  $$
  and hence
  $
  \sum_{n=1}^\infty H^2(\frac1{\gamma(B_n)}\,\gamma\restriction B_n, \mu_n)<   +\infty.
  $
  Therefore, by  \cite[Theorems 3.9, 3.6]{Hi},    $\gamma^{\boldsymbol B}\sim\bigotimes_{n=1}^\infty\mu_n$.
  Hence $\boldsymbol T$ is an IDPFT, as claimed.
  We now deduce from Proposition~\ref{pro2.9} that $\boldsymbol T$ is either sharply weak mixing or totally dissipative.
  Therefore, to complete the proof of (i)$\Rightarrow$(ii), it would be enough to show that $\boldsymbol T$
  is conservative.
  Unfortunately, we can not prove this fact.
However we observe  that if one replaces $\boldsymbol B$ with an arbitrary infinite subsequence 
  $\boldsymbol B'$ and associates a  $G$-action  $\boldsymbol T'$ with $\boldsymbol B'$ in the same way as we associated $\boldsymbol T$ with $\boldsymbol B$ then  $\boldsymbol T'$ possesses the same properties that we established for $\boldsymbol T$: it is free, 0-type, IDPFT and  it admits an exhausting $\boldsymbol T'$-F{\o}lner sequence.
  Therefore, to complete the proof of (i)$\Rightarrow$(ii), it suffices (in view of Proposition~\ref{pro2.9} to 
select a
subsequence $\boldsymbol B'$ of $\boldsymbol B$ such that  the corresponding $G$-action $\boldsymbol T'$  is conservative.

  Let $(\mathcal P_n)_{n=1}^\infty$ be a  sequence of finite  partitions of $Y$ into Borel subsets of equal measure
  such that the sequence $(\mathcal P_n)_{n=1}^\infty$ approximates $(\frak C,\gamma)$.   
  Fix a sequence $(K_n)_{n=1}^\infty$
 of compact subsets of $G$  such that $K_1\subset K_2\subset\cdots$ and $\bigcup_{n=1}^\infty K_n=G$.
 Since the $G$-action $(T_g^{\otimes n})_{g\in G}$ on $Y^n$ preserves the probability measure 
 $\gamma^{\otimes n}$, this action is conservative.
 Hence for each subset $A\subset Y^n$, there is $g\in G\setminus K_n$ such that
 $\gamma^{\otimes n}(T_g^{\otimes n}A\cap A)>0$.
 Moreover, if we fix a countable dense subgroup $G'$ of $G$ then we can additionally claim that
 $g\in G'$.
 Applying a standard exhaustion argument, we  construct a sequence $(A_m)_{m=1}^\infty$
 of mutually disjoint subsets of $A$ and a sequence $(g_m)_{m=1}^\infty$ of elements of $G'$ such that
 $A=\bigsqcup_{m=1}^\infty A_m=\bigsqcup_{m=1}^\infty T_{g_m}^{\otimes n}A_m$ and $g_m\in G'\setminus K_n$ for each $m\in\Bbb N$.
 Hence there is $M>0$ such that  
 \begin{gather*}
 \bigg(\bigsqcup_{m=1}^MA_n\bigg)\cup\bigg(\bigsqcup_{m=1}^MT_{g_m}^{\otimes n}A_n\bigg)\subset A\quad\text{  and}\\
  \gamma^{\otimes n}\bigg(\bigsqcup_{m=1}^MA_m\bigg)=\gamma^{\otimes n}\bigg(\bigsqcup_{m=1}^MT_{g_m}^{\otimes n}A_m\bigg)>0.5\gamma^{\otimes n}(A).
  \end{gather*}
  Then for every $n>0$, we apply this argument to each atom $P$ of the partition $(\mathcal P_n)^{\otimes n}$ of $Y^n$ to determine
 a finite subset $F_n\subset G\setminus K_n$ and  a family of measured subsets $(P_f)_{f\in F_n}$ of $P$ satisfying the following conditions:
 \begin{itemize}
 \item
  $P_f\cap P_h=\emptyset$ and $(T_f)^{\otimes n}P_f\cap (T_h)^{\otimes n}P_h=\emptyset$ if $f\ne h$,
  \item
  $\bigsqcup_{f\in F_n}(T_f)^{\otimes n}P_f\subset P$ and
  \item 
  $\gamma^{\otimes n}(\bigsqcup_{f\in F_n}P_f)>0.5\gamma^{\otimes n}(P)$.
  \end{itemize}
  The set $F_n$ is common for all atoms of $(\mathcal P_n)^{\otimes n}$ as we admit that some
  subsets $P_f$ can be of zero measure.
  We now select a sequence $k_1<k_2<\dots$ in such a way that
  $$
  \max_{f\in F_{k_n}}\prod_{m>n}\frac{\gamma(T_fB_{k_m}\cap B_{k_m})}{\gamma(B_{k_m})}>0.5
  \quad\text{ for each $ n$.}
  $$
  This is possible because $\boldsymbol B$ is $T$-asymptotically invariant.
  For simplicity of notation, we let $F_n':=F_{k_n}$, $B'_n:=B_{k_n}$,   $\boldsymbol B'=(B_n')_{n=1}^\infty$,
  $$
 (\boldsymbol B')^n:=Y^n\times B_{n+1}'\times B_{n+2}'\cdots,$$
   $\mathcal P_n':=\mathcal P_{k_n}$ and  $P_f':=P_f\in\mathcal P_n'$ if $f\in F_n'$.
  Then  $(\mathcal P_n')_{n=1}^\infty$ 
is a  sequence of finite  partitions of $Y$ into Borel subsets of equal measure and 
 the sequence $(\mathcal P_n')_{n=1}^\infty$ approximates $(\frak C,\gamma)$.
    We set  
   $$
   \mathcal P_n^*:=\{P'\times B_{n+1}'\times B_{n+2}'\times\cdots\mid P'\in(\mathcal P_n')^{\otimes n}\}.
   $$
   Then $\mathcal P_n^*$ is a finite partition of $(\boldsymbol B')^n$ into subsets of equal measure for each
   $n\in\Bbb N$.
   Moreover, 
   the sequence $( \mathcal P_n^*)_{n=1}^\infty$ approximates
   $(\frak B,\gamma^{\boldsymbol B'})$.

  Take  an atom $ P^*\in\mathcal P_n^*$ for some $n\in\Bbb N$.
  Then $P^*=P'\times  B_{n+1}'\times B_{n+2}'\times\cdots$ for an atom
   $P'\in(\mathcal P_n')^{\otimes n}$.
  We now set  for each $f\in F_n'$,
 $$
  P_f^*:=P_f'\times  (B_{n+1}'\cap T_f^{-1}B_{n+1}')\times (B_{n+2}'\cap T_f^{-1}B_{n+2}')\times\cdots.
  $$
  Then 
    $( P_f^*
    )_{f\in F_n'}$ are mutually disjoint Borel subsets of $ P^*$ and  $(\boldsymbol T'_fP_f^*)_{f\in F_n'  }$ are also
    mutually disjoint Borel subsets of $P^*$.
    Moreover,
    \begin{equation}\label{3-2}
  \begin{aligned}
 \gamma^{\boldsymbol B'}\bigg(\bigsqcup_{f\in F_n'}P_f^*\bigg)
 &=\sum_{f\in F_n'}\frac{\gamma^{\otimes n}(P_f')}{\prod_{j=1}^n\gamma(B_j')}\prod_{k>n}\frac{\gamma(T_fB_k'\cap B_k')}{\gamma(B_k')}\\
 &>\frac12\sum_{f\in F_n'}\frac{\gamma^{\otimes n}(P_f')}{\prod_{j=1}^n\gamma(B_j')}\\
 &=
\frac{\gamma^{\otimes n}(\bigsqcup_{f\in F_n'}P_f')}{2\prod_{j=1}^n\gamma(B_j')}\\
 &>
 \frac{ \gamma^{\boldsymbol B'}(P^*)}4.
 \end{aligned}
 \end{equation}
It follows from this and Proposition~\ref{pro3.1}(i) that $\boldsymbol T'$ is conservative, as desired.

To prove (i)$\Rightarrow$(iii), we first note that there exists a mixing measure preserving free $G$-action $T$ on a standard probability space $(Y,\frak C,\gamma)$ that admits, for each $m\in\Bbb N$, a 
$T$-asymptotically invariant sequence  $(A_{n,m})_{n=1}^\infty$ with  $\gamma(A_{n,m})=2^{-m}$ for all $n,m\in\Bbb N$.
\footnote{Indeed, let a  mixing action $S=(S_g)_{g\in G}$ on a standard probability space $(Z,\kappa)$ have an asymptotically invariant sequence  $(A_n)_{n=1}^\infty$ with $\kappa(A_n)=\frac12$ as in
\cite[Theorem~2.2.2]{Ch--Va}.
Consider an infinite product $G$-action $R=((S_g)^{\otimes\Bbb N})_{g\in G}$.
Of course, $R$ is  mixing.
For each $m>1$,  let $A_{n,m}:=(A_n)^m\times Y\times Y\times\cdots\subset Y^{\Bbb N}$.
The sequence $(A_{n,m})_{n=1}^\infty$ is $R$-asymptotically invariant and $\kappa^{\otimes\Bbb N}(A_{n,m})=2^{-m}$.
}
Utilizing this action $T$ and the sequence $\boldsymbol B:=(A_{n,1})_{n=1}^\infty$ we 
construct
the dynamical system  $(X,\frak B,\gamma^{\boldsymbol B'},\boldsymbol T')$ as above.
The  standard diagonalization argument applied to $(A_{n,m})_{n=1}^\infty$ yields a  $T$-F{\o}lner
sequence $(A_n)_{n=1}^\infty$
  such that
$\gamma (A_n)=2^{-n}$ for each $n\in\Bbb N$.
We now set 
$
\boldsymbol A^n:=Y^{n-1}\times A_n\times B'_{n+1}\times B'_{n+2}\times\cdots\subset X.
$
Then  
\begin{align*}
\gamma^{\boldsymbol B'}(\boldsymbol A^n)&=\frac{\gamma(A_n)}{\prod_{j=1}^n\gamma(B_j')}=1\qquad\text{ and}\\
  \gamma^{ \boldsymbol B'}(\boldsymbol T'_g\boldsymbol A^n\cap\boldsymbol A^n)
  &=\Bigg(\,{\prod_{j=1}^{n-1}\gamma(B_j')^{-1}}\Bigg)\frac{\gamma (A_n\cap T_g A_n)}{\gamma(B_n')}
 \prod_{j>n}\frac{\gamma(B_j'\cap T_gB_j')}{\gamma(B_j')}\\
 &=\frac{\gamma (A_n\cap T_g A_n)}{\gamma(A_n)}
 \prod_{j>n}\frac{\gamma(B_j'\cap T_gB_j')}{\gamma(B_j')}.
 \end{align*}
  Since $\sup_{g\in K} \prod_{j>n}\frac{\gamma(B_j'\cap T_gB_j')}{\gamma(B_j')}\to 1$ for each compact subset $K\subset G$ in view of~(\ref{3-1}),
 we obtain that $\sup_{g\in K}\gamma^{\boldsymbol B}(\boldsymbol T'_g\boldsymbol A^n\cap\boldsymbol A^n)\to 1$ as $n\to\infty$, as desired.
 \end{proof*}

\begin{proof*}{Proof of Corollary C} The ``if'' part follows from \cite[Theorem~2.2.2]{Ch--Va}.
 We prove the ``only if'' part. 
Let $G$ has the Haagerup property.
By Theorem~B,
there exists a  sharply weak mixing  0-type  measure preserving free  $G$-action $T$ on an infinite $\sigma$-finite standard measure space $(X,\frak B,\mu)$ admitting a  $T$-F{\o}lner sequence $(A_n)_{n=1}^\infty$ such that $\mu(A_n)=1$ for all $n\in\Bbb N$.
Denote by $(X^*,\frak B^*,\mu^*, T^*)$ the Poisson suspension of 
$(X,\frak B,\mu, T)$. 
 For each $n\in\Bbb N$, we set $[A_n]_0:=\{\omega\in X^*\mid \omega(A_n)=0\}$.
 Then 
 $$
\mu^*([A_n]_0)=e^{-{\mu(A_n)}}=e^{-1}
  $$
  and
 \begin{align*}
   \mu^*(T_g^*[A_n]_0\cap [A_n]_0)&=  \mu^*([T_gA_n\cup A_n]_0)\\
 &=e^{-\mu(\boldsymbol T_gA_n\cup A_n)}.
 \end{align*}
 Since for each compact subset $K\subset G$, 
 $$
 \sup_{g\in K}|\mu( T_g A_n\cup A_n)-\mu( A_n)|\to 0,
 $$
  we obtain that 
 $$
 \sup_{g\in K}|  \mu^*(T_g^*[A_n]_0\cap [A_n]_0)-\mu^*([A_n]_0)
|\to 0
 $$
  as $n\to\infty$.
 Thus, the sequence $([ A_n]_0)_{n=1}^\infty$ is  nontrivial and $ T^*$-asymptotically invariant.
  Thus $ T^*$ is not strongly ergodic.
  Since $ T$ is of 0-type, it follows from Fact~1.12(iii) that $U_{T^*}(g)\restriction L^2_0\big(X^*,\mu^*\big)\to 0$ weakly as $g\to\infty$.
  Hence 
  $ T^*$ is mixing.
  \end{proof*}

\section{Pairs of groups with Kazhdan property~(T)}\label{sec4}

 In this section we prove Theorem~D  and Corollary~E.
 Prior to this we note that if  $S$ is a weakly mixing measure preserving $G$-action on an infinite  $\sigma$-finite standard measure space $(Y,\frak F,\nu)$
 then the Koopman representation $U_S$ is weakly mixing.
 Indeed, suppose that $U_S$ contains a finite dimensional subspace.
 Then there is a unitary representation $V$ of $G$ in  a finite dimensional Hilbert space $\mathcal H$ and  a nontrivial mapping $F:Y\to\mathcal H$ such $F(S_gy)=V(g)F(y)$ for each $g\in G$ at a.e. $y\in Y$.
By \cite[Theorem~1.1]{GlWe2},  $F$ is constant a.e. 
 Since $\nu$ is infinite and
 the mapping $Y\ni y\mapsto\langle F(y),h\rangle$ belongs to $L^2(Y,\nu)$ for each $h\in\mathcal H$, it follows that $F=0$.
 Therefore $U_S$ is weakly mixing, as claimed.
 It now follows from Fact~\ref{fac2.12}(v) that $S$ is of weak 0-type.
 Therefore Theorem~D implies \cite[Theorem~1.5]{Jo2}.

 \begin{proof*}{Proof of Theorem~D}
 (i) If there is an $S$-F{\o}lner sequence then \ref{ast}\footnote{This is a formula on the first page of the paper.} holds for the  Koopman representation $U_S$ of $G$.
 Hence the restriction of $U_S$ to $H$ has a nontrivial invariant vector.
 Hence the action $S\restriction H$ has an invariant subset of positive finite measure.
 This contradicts to the condition of (i).
 
 (ii) Let $H\subset G$ do not have  property (T).
 The beginning of our argument is a slight modification of  the proof of \cite[Theorem~1.5]{Jo2}.
 There exists a conditionally negative definite function $\psi:G\to\Bbb R_+$ which is unbounded on $H$ \cite[Theorem~1.2(a4')]{Jo1}.
 By the Schoenberg theorem,  for each $t>0$, the function $\phi_t:=e^{-t^{-1}\psi}:G\to\Bbb R_+$ is positive definite. 
 Hence
the GNS-construction yields a triplet $(V_t,\mathcal H_t,\xi_t)$, consisting of a separable Hilbert space $\mathcal H_t$, a unitary representation $V_t$ of $G$ in $\mathcal H_t$ and a $V_t$-cyclic unit vector $\xi_t\in\mathcal H_t$ such that $\langle V_t(g)\xi_t,\xi_t
\rangle=\phi_t(g)$ for each $g\in G$. 
Since $\psi$ is unbounded on $H$, there is a sequence $(h_k)_{k=1}^\infty$ of elements in $H$ such that
$\psi(h_k)\to +\infty$.
As was shown in the proof of \cite[Lemma~2.1]{Jo1}, $\psi(g_1h_kg_2)\to +\infty$ as $k\to\infty$ for all $g_1,g_2\in G$.
Since 
$\xi_t$ is $V_t(G)$-cyclic,  it follows that $V_t(h_k)\to 0$ weakly as $k\to\infty$ for each $t>0$.
Thus,
the restriction of $V_t$ to $H$ is weakly mixing
by Fact~\ref{fac2.12}(v). 
Since $\phi_t$ takes only real values, $V_t$ is the  complexification of an orthogonal
representation of $G$. 
We will now argue as in the case (a) of the proof  of the main result from  \cite{CoWe}.
Let 
$$
\mathcal H:=\bigoplus_{n=1}^\infty \mathcal H_n\quad\text{ and }\quad V:=\bigoplus_{n=1}^\infty V_{n}.
$$
Of course, the restriction of $V$ to $H$ is  weakly mixing.
Moreover, $V$ is also the  complexification of an orthogonal
representation of $G$. 
Denote by $T=(T_g)_{g\in G}$ the corresponding Gaussian measure preserving $G$-action
on a standard probability space $(Y,\frak C,\gamma)$.
 By Fact~\ref{fac2.12}(iv), the associated Koopman representation
$U_T$  of $G$ is unitarily equivalent to   $\exp V$. 
Since $\exp (V\restriction H)=(\exp V)\restriction H$,
it follows from Fact~\ref{fac2.12}(ii)  that the unitary representation
$(U_T(h))_{h\in H}$
is weakly mixing on $L^2_0(Y,\gamma)$.
Hence  $(T_h)_{h\in H}$ is weakly mixing by Fact~\ref{fac2.12}(i).
We recall that $\mathcal H$ is a subspace of $\mathcal F(\mathcal H)$ which is canonically isomorphic to
$L^2(Y,\gamma)$.
Since $\mathcal H_n$ is a subspace of $\mathcal H$, we obtain that $\xi_n\in L^2(Y,\gamma)$, $\xi_n$ is a centered Gaussian variable on $Y$ and 
$$
\langle U_T(g)\xi_n,\xi_n\rangle=\langle V_n(g)\xi_n,\xi_n
\rangle=\phi_n(g)=e^{-\psi(g)/n}\qquad\text{for each }g\in G.
$$
 Let $B_n:=\{y\in Y\mid \xi_n(y)>0\}$.
 Denote by $\vartheta_n$ the distribution of  the random variable $\xi_n$.
 Then $\gamma(B_n)=\vartheta_n((0,+\infty))=0.5$ by the symmetry of $\vartheta_n$.
As in the proof of the main result from \cite{CoWe}, we obtain that
 $$
 \gamma(T_gB_n\triangle B_n)=
 \frac{\arccos  \langle U_T(g)\xi_n,\xi_n\rangle       }{\pi}= \frac{\arccos e^{-\psi(g)/n}}{\pi}.
  $$
This yields that  the sequence $\boldsymbol B:=(B_n)_{n=1}^\infty$ is $T$-asymptotically invariant. 
Moreover,
without loss of generality we may (and will) assume that (\ref{3-1}) holds for $\boldsymbol B$.

Consider now the dynamical system $(X,\frak B,\gamma^{\boldsymbol B},\boldsymbol T)$ as in the proof of Theorem~B.
Then $\gamma^{\boldsymbol B}(X)=\infty$ and $\boldsymbol T$ preserves $\gamma^{\boldsymbol B}$
by Proposition~\ref{pro2.6}.
As  was shown in the proof of Theorem~B,   the sequence $(\boldsymbol B^n)_{n\in \Bbb N}$ (defined there) is $\boldsymbol T$-F{\o}lner and exhausting. 

Reasoning as in the proof of Theorem~B,
we have only to choose  a  subsequence $\boldsymbol B'$ in $\boldsymbol B$ in such a way that the restriction of  $\boldsymbol T'$ (we recall  again that $\boldsymbol T'$ is determined by $\boldsymbol B'$) to $H$
 is weakly mixing.
Let $(\mathcal P_n)_{n=1}^\infty$  denote the same sequence of finite partitions as in the proof of Theorem~B.
Since $(T_h)_{h\in H}$ is weakly mixing, the $H$-action $(T_h^{\otimes{2n}})_{h\in H}$ on the probability space $(X^{2n},\gamma^{\otimes 2n})$ is ergodic.
Hence for every $n>0$, there is a subset $F_n\subset H$ such that for every two atoms $P,Q$
of the partition $(\mathcal P_n)^{\otimes 2n}$ of $B_n\times B_n$, there is a family of measurable subsets $(P_f)_{f\in F_n}$ of $P$ such that
\begin{itemize}
 \item
  $P_f\cap P_h=\emptyset$ and $(T_f)^{\otimes 2n}P_f\cap (T_f)^{\otimes 2n}P_f=\emptyset$ if $f\ne h$,
  \item
   $\bigsqcup_{f\in F}(T_f)^{\otimes 2n}P_f\subset Q$ and
  \item
   $\gamma^{\otimes 2n}(\bigsqcup_{f\in F_n}P_f)>0.5\gamma^{\otimes2 n}(P)$.
\end{itemize}
We now select a sequence $k_1<k_2<\dots$ of integers in such a way that
  $$
  \max_{f\in F_{k_n}}\prod_{m>n}\frac{\gamma(T_fB_{k_m}\cap B_{k_m})}{\gamma(B_{k_m})}>0.5
  \quad\text{ for each $ n$.}
  $$
As in the proof of Theorem~B, we let $F_n':=F_{k_n}$, $B'_n:=B_{k_n}$, $\boldsymbol B'=(B_n')_{n=1}^\infty$, $(\boldsymbol B')^n:=Y^n\times B_{n+1}'\times B_{n+2}'\cdots$, $\mathcal P_n':=\mathcal P_{k_n}$ and  $P_f':=P_f\in\mathcal P_n'$ if $f\in F_n'$.
We also set 
\begin{align*}
 P^* &:=P'\times (B'_{n+1}\times B_{n+2}'\times\cdots)^{\times 2}\subset(\boldsymbol B')^n\times(\boldsymbol B')^n,\\
 Q^*&:=Q'\times (B_{n+1}'\times B_{n+2}'\times\cdots)^{\times 2}\subset(\boldsymbol B')^n\times(\boldsymbol B')^n
 \qquad \text{and}\\ 
 P_f^*&:=P_f'\times\big((B_{n+1}'\cap T_f^{-1}B_{n+1}')\times (B_{n+2}'\cap T_f^{-1}B_{n+2}')\times\cdots\big)^{\times 2}
 \end{align*}
 for all atoms $P',Q'\in (\mathcal P_n')^{\otimes 2n}$.
 Then $\mathcal P_n^*:=(P^*)_{P'\in(\mathcal P_n')^{\otimes 2n}}$ is a finite partition of $(\boldsymbol B')^n\times(\boldsymbol B')^n$ into subsets of equal measure
and the sequence $( \mathcal P_n')_{n=1}^\infty$ approximates
$(\frak B\otimes\frak B,\gamma^{\boldsymbol B'}\otimes \gamma^{\boldsymbol B'})$.
 We also note that
$(P^*_f)_{f\in F_n}$ are  mutually disjoint subsets of $P^*$
and $((\boldsymbol T_f\times\boldsymbol T_f)P^*_f)_{f\in F_n}$ 
are  mutually disjoint subsets of $Q^*$.
Arguing as in \ref{3-2}, we obtain that
$$
(\gamma^{\boldsymbol B'}\otimes\gamma^{\boldsymbol B'})\bigg(\bigsqcup_{f\in F_n}P_f^*\bigg)>
\frac18(\gamma^{\boldsymbol B'}\otimes\gamma^{\boldsymbol B'})(P^*).
$$
It follows from that and Proposition~\ref{pro3.1}(ii) that the $H$-action $(\boldsymbol T'_g\times\boldsymbol T'_g)_{g\in H}$ is ergodic.
Hence  $\boldsymbol T'\restriction H$ is weakly mixing by \cite[Theorem~1.1]{GlWe2}.
  \end{proof*}

 Corollary~E follows from Theorem~D in the same way as Corollary~C follows from Theorem~B.

  \section*{Acknowledgements}
  
I am grateful to the anonymous referee for the careful reading the manuscript and his (her) numerous useful remarks.

\end{document}